\documentclass[preprint,3p,11pt, times, authoryear]{elsarticle}
 
\usepackage{booktabs}
\usepackage{float}

\usepackage{amsmath,amssymb,amsfonts,mathrsfs,hyperref,microtype, amsthm}

\usepackage[T1]{fontenc}

\usepackage{lmodern} 

\usepackage{booktabs}
\usepackage{float}
\usepackage{amsmath}
\usepackage{amsthm}
\usepackage{amsmath,amssymb,amsfonts,mathrsfs,hyperref,microtype, amsthm}
\usepackage{graphicx}
\usepackage{graphicx,mathabx}
\usepackage{enumerate}
\usepackage{amsmath}
\usepackage{txfonts}
\usepackage{soul}
\usepackage{xcolor,graphicx}
\AtBeginDocument{\mathcode`v=\varv}

\usepackage{subfig}
\graphicspath{ {images/} }
\numberwithin{equation}{section}

\newcommand{\BS}{\boldsymbol}

\newcommand{\Rmnum}[1]{\expandafter\@slowromancap\romannumeral #1@}
\newcommand{\trp}{{\sf T}}

\journal{Journal of Statistical Planning and Inference}

\begin{document}

\begin{frontmatter}
\author[1]{Osama Idais\corref{cor1}}
\ead{osamaidais@gmail.com}
\cortext[cor1]{Corresponding author}
\address[1]{Independent Researcher, Mannheim, Germany}

\title{On the robustness of locally \(D\)-optimal designs under parameter uncertainty in single-factor generalized linear models}

\begin{abstract}
Local and Bayesian approaches to optimal experimental
design specify information about the regression parameters
through a nominal value and a prior distribution,
respectively. A practical question is whether accounting
for parameter uncertainty justifies replacing a locally
optimal design. We investigate when locally $D$-optimal
designs remain Bayesian $D$-optimal in single-factor
Poisson, gamma and logistic regression on the unit
interval. Analytical arguments and numerical comparisons
under normal and uniform slope priors examine how prior
form and dispersion affect the optimal support and the
Bayesian $D$-efficiency of retaining a local design.
\end{abstract}


\end{frontmatter}

\section{Introduction}
\label{sec1}

Generalized linear models (GLMs), introduced by \citet{10.2307/2344614}, include Poisson regression for counts,
gamma regression for positive continuous responses and logistic regression for binary responses. Optimal experimental design
selects factor settings and observation allocations to improve statistical precision. The $D$-optimality criterion maximizes the determinant of the Fisher information matrix, and equivalently minimizes the asymptotic variance-covariance matrix of the parameter estimators, see \cite{pukelsheim2006optimal}. In GLMs, this information matrix generally depends on unknown parameters, so the design optimization problem requires prior information about their values. The local approach was introduced by \citet{chernoff1953} to address this dependence by specifying a best guess parameter value. The resulting locally optimal design may change when this value changes.  Analytical results for locally optimal designs in GLMs are available (see \cite{10.2307/2346142, 10.2307/24308852, GAFFKE2019,  idais2021local}). 

Alternatively, the Bayesian approach addresses this parameter dependency by incorporating a prior distribution and maximizing the expected log-determinant of the information matrix \cite{ChalonerVerdinelli1995,atkinson2015designs}. Bayesian design for logistic regression was studied by \citet{chaloner1989optimal}, while \citet{van2026analytic} derived analytical expressions for Bayesian $D$-optimality criteria in Poisson regression. A related approach, proposed by \citet{rios2025bootstrap}, aggregates local GLM designs obtained from bootstrap samples of pilot data.

Both the local and Bayesian approaches require prior information about the regression parameters. The local
approach uses a nominal value, while the Bayesian approach uses a prior distribution. It is therefore 
natural to ask whether accounting for uncertainty through a prior distribution requires a different design and if the design changes, a further question is how much efficiency is gained by replacing the local design.

In this paper, we examine these questions for single-factor GLMs including Poisson, gamma and logistic regression
on the unit interval. We study when local $D$-optimality is preserved and assess the efficiency of retaining the local design under different prior specifications.
In particular, normal and uniform slope priors are used to examine the effects of prior form and dispersion on the optimal support. The comparison is conducted within the saturated class and over all designs. 

The paper is organized as follows.
Section~\ref{sec2} introduces the model framework and
the optimality criteria. Section~\ref{sec3} presents
analytical results and numerical comparisons within
the saturated class and over all designs.
Section~\ref{sec5} discusses the main findings.
All numerical computations were performed in Python.

\section{Preliminaries}
\label{sec2}

Consider independent responses $Y$  following a single-factor
GLM with linear predictor $\eta(x,\BS{\beta})=\beta_0+\beta_1x
=\mathbf{f}(x)^\trp\BS{\beta}$, where  $x$ comes for the unit interval experimental region $[0,1]$, the regression functional-vector  $\mathbf{f}(x)=(1,x)^\trp$  and the regression parameter vector $\BS{\beta}=(\beta_0,\beta_1)^\trp$ belongs to the
admissible parameter space $\mathcal{B}$. Denote the one-to-one link function $g$, then the mean is $\mu(x,\BS{\beta})=g^{-1}\{\eta(x,\BS{\beta})\}$,
and the intensity function is
\[
u(x,\BS{\beta})
=
\frac{1}{V\{\mu(x,\BS{\beta})\}}
\left(\frac{d\mu(x,\BS{\beta})}{d\eta}\right)^2,
\]
where  $V$ is the mean-variance function. Any common dispersion factor is omitted,
as it does not affect the optimal designs. For more details on GLMs, see \citet[Section 2.2.2]{mccullagh1989generalized}.

An approximate design $\xi$ is a probability measure
on $[0,1]$ (see \citet[p. 15]{silvey1980optimal}). For a design with distinct support points
$x_1,\ldots,x_m$ and weights $\omega_i>0$ satisfying
$\sum_{i=1}^m\omega_i=1$, its normalized information
matrix is
\begin{equation}
\mathbf{M}(\xi;\BS{\beta})
=
\int_0^1 u(x,\BS{\beta})
\mathbf{f}(x)\mathbf{f}(x)^\trp\,\xi(dx)
=
\sum_{i=1}^m\omega_i u(x_i,\BS{\beta})
\mathbf{f}(x_i)\mathbf{f}(x_i)^\trp.
\label{infoM}
\end{equation}
A design is saturated if it has two support points,
the number of regression parameters.

Throughout, a prior $\pi$ is called admissible if it is
proper, satisfies $\pi(\mathcal{B})=1$, and meets the
following integrability conditions. Therefore, we assume 
 we assume
$E_\pi(|\beta_0|+|\beta_1|)<\infty$ for the  Poisson and logistic regression. Also, for gamma regression under identity and reciprocal links,
we assume
$E_\pi\{|\log\beta_0|+|\log(\beta_0+\beta_1)|\}<\infty$.
These conditions ensure that the Bayesian criterion
is finite. For an admissible prior $\pi$, the Bayesian
$D$-criterion is
\begin{equation}
\phi_D(\xi;\pi)
=
E_\pi\!\left[
\log\det\mathbf{M}(\xi;\BS{\beta})
\right],
\label{eq:bayes-D-criterion}
\end{equation}
and a design maximizing this criterion is Bayesian
$D$-optimal. Local $D$-optimality at a nominal value
$\widetilde{\BS{\beta}}$ is recovered by taking
$\pi=\delta_{\widetilde{\BS{\beta}}}$. That is $E_\pi\!\left[
\log\det\mathbf{M}(\xi;\BS{\beta})
\right]=\log\det\mathbf{M}(\xi;\widetilde{\BS{\beta}})$. 
For any fixed two-point support, both criteria give
equal weights, since the determinant is proportional
to $\omega_1\omega_2$.

Define the Bayesian sensitivity function by
\begin{equation}
d_\pi(x,\xi)
=
E_\pi\!\left[
u(x,\BS{\beta})\mathbf{f}(x)^\trp
\mathbf{M}(\xi;\BS{\beta})^{-1}\mathbf{f}(x)
\right].
\label{D-sen}
\end{equation}
By the Bayesian equivalence theorem (Kiefer-Wolfowitz Equivalence Theorem, see  \citet{silvey1980optimal}, \citet{atkinson2015designs}), for an admissible
prior $\pi$, a design $\xi^*$ with nonsingular information matrix $\pi$-almost surely is Bayesian
$D$-optimal if and only if
\begin{equation}
d_\pi(x,\xi^*)\leq 2
\qquad\mbox{for all }x\in[0,1].
\label{eq:bayes-D-equivalence}
\end{equation}
Equality holds at every support point of $\xi^*$.

\section{Main results}
\label{sec3}
 In this section, we examine when a locally D-optimal design on $[0,1]$ remains optimal under parameter uncertainty, represented by  a prior distribution $\pi$ defined on admissible parameter space  $\mathcal B \subseteq \mathbb{R}^2$. A minimally supported design $\xi_{\omega}$ assigns weights $\omega$ and $1-\omega$ to two design points $x_1$ and $x_2$, $0\le x_1<x_2\le1$. Its information matrix satisfies $\det\mathbf M(\xi_{\omega},\BS{\beta})=\omega(1-\omega)(x_2-x_1)^2u(x_1,\BS{\beta})u(x_2,\BS{\beta})$. As is well known, the optimal weights of the optimal designs are $\omega=1/2$ under both the local and Bayesian $D$-criteria.  Denoting the equally weighted design by $\xi_{1/2}$, the Bayesian criterion reduces to
\begin{align}
\phi_D(\xi_{1/2};\pi)
&=
\log\!\Bigl(\tfrac14(x_2-x_1)^2\Bigr)
+\int_{\mathcal B}\log\!\bigl(u(x_1,\BS{\beta})u(x_2,\BS{\beta})\bigr)\pi(d\BS{\beta}) \label{det2}\\
&=
\log\!\Bigl(\tfrac14(x_2-x_1)^2\Bigr)
+E_\pi\!\left[\log u(x_1,\BS{\beta})+\log u(x_2,\BS{\beta})\right]. \notag
\end{align}
Thus, within the class of minimally supported designs, the Bayesian $D$-optimality problem reduces to the choice of $x_1$ and $x_2$. However, optimality within this class does not necessarily imply Bayesian $D$-optimality over all possible designs on $[0,1]$, as this relies on the form of the intensity function itself.

For illustration, we consider the gamma, Poisson, and logistic models. Under gamma model, both the identity link, $\mu(x,\BS{\beta})=\eta$, and the reciprocal link, $\mu(x,\BS{\beta})=1/\eta$, provide the same intensity $u(x,\BS{\beta})=1/\eta(x,\BS{\beta})^2$, and requires the admissible parameter space $\mathcal B=\{\BS{\beta}\in\mathbb R^2:\beta_0>0,\ \beta_0+\beta_1>0\}$, which ensures positivity of the mean $\mu=g^{-1}(\eta)$ on the design region $[0,1]$. The results of \cite{GAFFKE2019} showed that the locally $D$-optimal design on such a compact interval $[0,1]$ assigns equal weights to its vertices $0$ and $1$ for any $\BS{\beta}\in\mathcal B$. In the Bayesian approach, the criterion at this design is
$\phi_D(\xi_{1/2}^*;\pi)=\log(1/4)-2E_\pi[\log\beta_0+\log(\beta_0+\beta_1)]$. Since the equally weighted endpoint design is locally $D$-optimal at every admissible parameter value, it remains Bayesian $D$-optimal under every admissible prior. Thus, the optimal design is independent of the prior and account for uncertainty in those parameters does not change the optimal design for these one-factor gamma models. 

For Poisson model with a log link, $u(x,\BS{\beta})=\mu(x,\BS{\beta})=\exp\{\eta\}$, $\BS{\beta}\in\mathbb R^2$. The intercept term $\beta_0$ is independent of the design optimizations. As shown by \cite{atkinson2015designs}, the Bayesian $D$-criterion for minimally supported designs coincides with the local criterion evaluated at the prior mean $E_\pi(\BS{\beta})$. For this one-factor model, denote prior mean of the slope by $m=E_\pi(\beta_1)$,  the $D$-criterion becomes $\phi_D(\xi_{1/2};\pi)=\log((1/4) (x_2-x_1)^2)+2E_\pi(\beta_0)+m(x_1+x_2)$. The maximization reduces to $\max_{0\le x_1<x_2\le1}\{2\log(x_2-x_1)+m(x_1+x_2)\}$. Thus, combining this with \cite[Remark~3]{10.2307/24308852}, the optimal support, within the class of  minimally supported designs, depends on the prior only through $m$ and is given by
\[
(x_1^*,x_2^*)
=
\begin{cases}
\left(0,-\dfrac{2}{m}\right), & m<-2,\\[1.2ex]
(0,1), & -2\le m\le2,\\[1.2ex]
\left(1-\dfrac{2}{m},1\right), & m>2.
\end{cases}
\]
In this situation,  when the nominal slope value equals the prior mean, i.e., $\beta_1=\widetilde\beta_1=m$, the design is locally $D$-optimal at $\widetilde\beta_1$.  
This design remains Bayesian $D$-optimal within the class of minimally supported, two-point, designs. For example, considering both priors $\beta_1\sim N(3,1/3)$,  and $\beta_1\sim U(2,4)$ yield $m=3$, and hence the optimal two-point design assigns equal weights to the support $\{1/3,1\}$. Let $\xi_{\mathrm{L}}^*$ be the locally
$D$-optimal design at a nominal slope
$\widetilde{\beta}_1$, and let $\pi$ satisfy
$E_\pi(|\beta_0|+|\beta_1|)<\infty$.
Writing $m=E_\pi(\beta_1)$, the design
$\xi_{\mathrm{L}}^*$ remains Bayesian $D$-optimal
among the saturated designs if and only if
\begin{equation}
\begin{cases}
|m|\leq 2,
& \text{if }|\widetilde{\beta}_1|\leq 2,\\
m=\widetilde{\beta}_1,
& \text{if }|\widetilde{\beta}_1|>2.
\end{cases}
\label{eq:poisson-preservation}
\end{equation}
Thus, matching the nominal slope to the prior mean
is sufficient for preservation, but is necessary
only when $|\widetilde{\beta}_1|>2$.
Optimality over all designs requires the additional
verification by the Bayesian equivalence theorem condition~\ref{eq:bayes-D-equivalence}.

For logistic model, $\mu(x,\BS{\beta})=\exp\{\eta\}/(1+\exp\{\eta\})$ and $u(x,\BS{\beta})=\mu(1-\mu)$, $\BS{\beta}\in\mathbb R^2$. Unlike the Poisson case, the expected log-intensities $E_\pi\!\left[\log u(x_1,\BS{\beta})+\log u(x_2,\BS{\beta})\right]$, in \eqref{det2}, generally depends on features of the joint prior beyond its means.  
Since
$|\partial\log u(x,\BS{\beta})/\partial x|\leq|\beta_1|$
and admissibility ensures $E_\pi|\beta_1|<\infty$,
differentiation of the expected log-intensities
in \eqref{det2} may be performed
under the prior expectation.
 By differentiating the criterion $\phi_D(\xi_{1/2};\pi)$ with respect to $x_1$ and $x_2$, one finds that an interior optimum $0<x_1^*<x_2^*<1$ satisfies this system of two equations
\begin{align}
-\frac{2}{x_2^*-x_1^*}
-
E_\pi\!\left[
\beta_1\tanh\!\left(\frac{\beta_0+\beta_1x_1^*}{2}\right)
\right]
&=0, \label{eq-log1}\\[1ex]
\frac{2}{x_2^*-x_1^*}
-
E_\pi\!\left[
\beta_1\tanh\!\left(\frac{\beta_0+\beta_1x_2^*}{2}\right)
\right]
&=0. \label{eq-log2}\
\end{align}

At $x_1^*=0$, the l.h.s. of equation (\ref{eq-log1}) must be nonpositive and at $x_2^*=1$, the l.h.s of equation (\ref{eq-log2}) must be nonnegative.
For each fixed $\BS{\beta}$, the second derivative 
$\partial^2\log u(x,\BS{\beta})/\partial x^2
=-\beta_1^2/
[2\cosh^2\{(\beta_0+\beta_1x)/2\}]\leq 0$.
Thus, the expected log-intensities are concave in
their respective support points. Together with the
concavity of $2\log(x_2-x_1)$, this implies that the
Bayesian criterion is concave on
$0\leq x_1<x_2\leq 1$.
This means that, the interior stationarity conditions,
together with the stated boundary inequalities (\ref{eq-log1}) and (\ref{eq-log2}),
are necessary and sufficient for optimality within
the class of two-point designs.  However, the locally optimal design at \(\widetilde{\boldsymbol\beta}=E_\pi(\boldsymbol\beta)\) need not be Bayesian optimal, even within the class of minimally supported designs.  We next provide further numerical illustrations on the design's robustness against parameter uncertainty.

\subsection{Numerical comparisons}
We illustrate when a locally $D$-optimal design remains unchanged under parameter uncertainty for the one-factor gamma, Poisson, and logistic models on the unite interval $[0,1]$. We first optimize within the class of saturated designs and then extend the comparison to all designs on $[0,1]$.  When the locally \(D\)-optimal design fails to hold as Bayesian optimal, we evaluate its efficiency to quantify the resulting information loss. 

Let the nominal parameter vector be given by $\widetilde{\BS{\beta}}=(3,3)^\trp$.  The equally-weighted locally $D$-optimal design has support $\{0,1\}$ for gamma model, $\{1/3,1\}$ for Poisson model, and approximately $\{0,0.6754\}$ for logistic model. These local designs are kept fixed throughout the comparison. 

To evaluate robustness, we utilize the six prior specifications listed in Table~\ref{tab:priors_onefactor}, including uniform, normal, and mixed uniform-normal distributions with varying degrees of variance. The components $\beta_0$ and $\beta_1$ are independent under the distributions listed in the table.  All six priors have mean vector $(3,3)^\trp$ of $(\beta_0, \beta_1)^\trp$. This allows us to examine whether changes in prior dispersion and distributional form affect the local design at this nominal value  $\widetilde{\BS{\beta}}=(3,3)^\trp$.

Such prior comparisons for gamma model with an inverse link were considered by \cite{Habib2017} on a design region with boundaries depend on the regression parameters.   In contrast, for gamma model in the present work,  the prior distributions are restricted and normalized on $\mathcal B=\{\BS{\beta}\in\mathbb R^2:\beta_0>0,\ \beta_0+\beta_1>0\}$. This restriction generally changes the prior means and introduces dependence between the components. This delivers no changes on robustness of the locally optimal design.

\begin{table}[H]
\centering
\caption{Prior specifications for the one-factor comparison. The components are independent before any restriction to the gamma admissible parameter space.}
\label{tab:priors_onefactor}
\begin{tabular}{@{}lll@{}}
\toprule
Prior & Class & Specification \\
\midrule
$\pi_1$ & Uniform & $\beta_0\sim U(2,4),\quad \beta_1\sim U(2,4)$ \\
$\pi_2$ & Normal & $\beta_0\sim N(3,1/3),\quad \beta_1\sim N(3,1/3)$ \\
$\pi_3$ & Mixed & $\beta_0\sim U(2,4),\quad \beta_1\sim N(3,1/3)$ \\
$\pi_4$ & Diffuse uniform & $\beta_0\sim U(-2,8),\quad \beta_1\sim U(-2,8)$ \\
$\pi_5$ & Diffuse normal & $\beta_0\sim N(3,100/12),\quad \beta_1\sim N(3,100/12)$ \\
$\pi_6$ & Diffuse mixed & $\beta_0\sim U(-2,8),\quad \beta_1\sim N(3,100/12)$ \\
\bottomrule
\end{tabular}
\end{table}

\subsubsection{Two-point Bayesian designs}

For saturated designs, $\omega_1=\omega_2=1/2$, and optimization reduces to choosing $0\le x_1<x_2\le1$ to maximize the Bayesian $D$-criteria $\phi_D(\xi_{1/2};\pi)$ given in (\ref{det2}). The gamma and Poisson solutions follow from the preceding analytical results. For logistic model, we optimize the expected criterion numerically using Gaussian quadrature. Table~\ref{tab:onefactor_results} reports the resulting support points.

\begin{table}[H]
\centering
\caption{Bayesian $D$-optimal support points within the saturated class. The local reference designs are evaluated at $\widetilde{\BS{\beta}}=(3,3)^\trp$.}
\label{tab:onefactor_results}
\renewcommand{\arraystretch}{1.15}
\begin{tabular}{@{}llll@{}}
\toprule
Prior & Gamma & Poisson & Logistic \\
\midrule
Local reference & $(0,1)$ & $(1/3,1)$ & $(0,0.67547)$ \\
\midrule
$\pi_1$ & $(0,1)$ & $(1/3,1)$ & $(0,0.67694)$ \\
$\pi_2$ & $(0,1)$ & $(1/3,1)$ & $(0,0.67697)$ \\
$\pi_3$ & $(0,1)$ & $(1/3,1)$ & $(0,0.67692)$ \\
$\pi_4$ & $(0,1)$ & $(1/3,1)$ & $(0,0.68402)$ \\
$\pi_5$ & $(0,1)$ & $(1/3,1)$ & $(0,0.68725)$ \\
$\pi_6$ & $(0,1)$ & $(1/3,1)$ & $(0,0.68072)$ \\
\bottomrule
\end{tabular}
\end{table}
As expected, the numerical solutions for the gamma and Poisson models agree with the analytical results. The local result for the gamma model has support  $\{0,1\}$ that remains unchanged under every admissible prior considered and is Bayesian $D$-optimal over all designs. 

The Poisson model allows all six priors to have mean slope $m=3$. Hence, the Bayesian optimal support within the class of saturated designs is $(1/3,1)$, which coincides with the local support at $\widetilde{\BS{\beta}}=(3,3)^\trp$. Within this class, under a fixed  mean slope, the optimal design is robust against changing the prior dispersion or distributional form. However, global optimality over all possible designs requires a separate verification by the equivalence theorem.

For logistic model, all six Bayesian solutions retain the lower support point $x_1^*=0$, but the upper support point differs from the local value $0.67547$. These differences are small, particularly under $\pi_1$, $\pi_2$, and $\pi_3$. Equality of prior means does not necessarily give the same optimal saturated design.

\subsubsection{Bayesian optimality over all designs}

We next consider Bayesian $D$-optimal designs over all designs on $[0,1]$, allowing more than two support points. Here, prior expectations were evaluated by Gaussian quadrature, using Gauss--Hermite rules for normal priors and Gauss--Legendre rules for uniform priors. For priors with independent uniform and normal components, the
corresponding rules were combined across the parameter coordinates. The computations used $50$ and $120$ nodes per coordinate, with $160$ nodes used for verification.

Numerical optimization for Tables~\ref{tab:onefactor_results}--\ref{tab:local_efficiency_joint} was performed using SLSQP in SciPy. For saturated logistic designs,
the two ordered support points were optimized with equal weights. For comparisons over all designs on
$[0,1]$, candidate designs supported on $\{0,x,1\}$ were optimized over the interior point and nonnegative
weights satisfying $\sum_i\omega_i=1$. Allowing zero weights includes two-point designs in this candidate
class. Multiple starting designs were used, and among the converged runs, the solution with the largest
Bayesian $D$-criterion value was retained. 

Optimality was assessed numerically using the Bayesian equivalence condition~\eqref{eq:bayes-D-equivalence}
and the sensitivity function~\eqref{D-sen},  $d_\pi(x,\xi^*)=E_\pi\!\left[u(x,\BS{\beta})\mathbf f(x)^\trp\mathbf M(\xi^*,\BS{\beta})^{-1}\mathbf f(x)\right]$ with $\mathbf f(x)=(1,x)^\trp$. For both Poisson and logistic regression, the maximum computed sensitivity $d_\pi(x,\xi^*)$ did not exceed $2+2\times10^{-7}$.

Under $\pi_1$, $\pi_2$ and $\pi_3$, these checks confirmed that the saturated Bayesian designs in Table~\ref{tab:onefactor_results} were also optimal  over all designs on $[0,1]$, within the stated numerical tolerance. Under the diffuse priors $\pi_4$, $\pi_5$ and $\pi_6$, the Bayesian optimal designs for both
models had three support points with unequal weights, as reported in Table~\ref{tab:onefactor_unrestricted}.

The designs under gamma models remains unchanged, supported by $0$ and $1$, under all six priors.

\begin{table}[H]
\centering
\caption{Bayesian $D$-optimal designs under the diffuse priors. Weights correspond to the support points in the order shown.}
\label{tab:onefactor_unrestricted}
\small
\renewcommand{\arraystretch}{1.15}
\setlength{\tabcolsep}{4pt}
\begin{tabular}{@{}llll@{}}
\toprule
Model & Prior & Support points & Weights \\
\midrule
Poisson & $\pi_4$ & $(0,\ 0.63949,\ 1)$ & $(0.25955,\ 0.28406,\ 0.45640)$ \\
        & $\pi_5$ & $(0,\ 0.61140,\ 1)$ & $(0.25023,\ 0.28842,\ 0.46135)$ \\
        & $\pi_6$ & $(0,\ 0.61140,\ 1)$ & $(0.25023,\ 0.28842,\ 0.46135)$ \\
\midrule
Logistic & $\pi_4$ & $(0,\ 0.38717,\ 1)$ & $(0.45680,\ 0.27053,\ 0.27267)$ \\
         & $\pi_5$ & $(0,\ 0.40057,\ 1)$ & $(0.45809,\ 0.26964,\ 0.27227)$ \\
         & $\pi_6$ & $(0,\ 0.40587,\ 1)$ & $(0.46056,\ 0.26981,\ 0.26962)$ \\
\bottomrule
\end{tabular}
\end{table}

For Poisson model, the local design at the prior mean remains optimal under the concentrated priors, but not under the diffuse priors. Although all six priors have the same  mean, the corresponding Bayesian optimal designs over all designs on $[0,1]$ need not coincide. The reduction to the prior mean therefore applies within the class of saturate designs and cannot generally be extended to all designs. The optimal design under priors $\pi_5$ and $\pi_6$ are equivalent because the slope priors are the same. 

For logistic model, the concentrated priors, $\pi_1$, $\pi_2$ and $\pi_3$, lead to small changes in the local support, whereas the diffuse priors change both the number of support points and their weights. To assess the practical importance of these changes, let $\xi_{\mathrm L}^*$ denote the locally D-optimal design at $\widetilde{\BS{\beta}}=(3,3)^\trp$ and $\xi_{\mathrm B}^*(\pi)$ the Bayesian optimal design under $\pi$. The Bayesian efficiency of retaining the local design is
\[
\operatorname{Eff}_{D,\pi}(\xi_{\mathrm L}^*)
=
\exp\!\left\{
\frac{\phi_D(\xi_{\mathrm L}^*;\pi)
-\phi_D(\xi_{\mathrm B}^*(\pi);\pi)}{2}
\right\},
\]
where
 $
 \phi_D(\xi_{\mathrm L}^*;\pi)
=E_\pi[\log\det M(\xi_L^*;\boldsymbol\beta)]$ 
and 
$
\phi_D(\xi_{\mathrm B}^*(\pi);\pi)
=E_\pi[\log\det M(\xi_B^*(\pi);\boldsymbol\beta)].
$
Both designs are evaluated under the same prior $\pi$. An efficiency of one means that the local design is also Bayesian $D$-optimal under the specified prior. 
A value close to one indicates that retaining the local design causes little loss, even when the Bayesian optimal design differs. 

The efficiency results showed that, for gamma model, this efficiency is exactly one under every prior considered. For Poisson model, it equals one under $\pi_1$, $\pi_2$, and $\pi_3$. For logistic model, it exceeds $99.99\%$ under these same three priors, so the small changes on the support have negligible impact on efficiency. Under the diffuse priors, however, the efficiency losses become more substantial, as shown in Table~\ref{tab:local_efficiency_joint}.

\begin{table}[H]
\centering
\caption{Bayesian $D$-efficiency, in $\%$, of the locally $D$-optimal design at $\widetilde{\BS{\beta}}=(3,3)^\trp$ relative to the Bayesian optimal design under the same prior.}
\label{tab:local_efficiency_joint}
\begin{tabular}{@{}lccc@{}}
\toprule
Prior & Gamma & Poisson & Logistic \\
\midrule
$\pi_4$ & 100 & 90.95 & 93.02 \\
$\pi_5$ & 100 & 91.89 & 93.77 \\
$\pi_6$ & 100 & 91.89 & 94.37 \\
\bottomrule
\end{tabular}
\end{table}

\subsubsection{Efficiency under slope uncertainty}

The previous results showed that diffuse priors can
reduce the efficiency of the local designs.
To examine the effect of slope uncertainty more directly,
we retain the same local reference designs and assign
a point-mass prior to $\beta_0=3$. We consider
$\beta_1\sim N(3,\sigma^2)$ and
$\beta_1\sim U(3-\sqrt{3}\sigma,3+\sqrt{3}\sigma)$,
which have the same mean and variance at each $\sigma$.
The case $\sigma=0$ denotes a point mass at 3.
The intercept is fixed only in the design prior;
both regression coefficients remain parameters to be
estimated, so the information matrix remains
two-dimensional.

For gamma regression, positivity of the mean throughout
$[0,1]$ requires $\beta_1>-3$.
Normal slope priors with $\sigma>0$ are therefore
conditioned on this region, with $\sigma$ denoting
the standard deviation before conditioning.
The uniform slope priors are already admissible
for $0\leq\sigma\leq3$.
The equally weighted endpoint design remains Bayesian
$D$-optimal, with efficiency exactly one.

For Poisson and logistic regression, candidate designs
supported on $\{0,x,1\}$ were optimized over the interior
point and nonnegative weights using SLSQP with multiple
starting designs. Their optimality over all designs on
$[0,1]$ was assessed numerically using the Bayesian
equivalence condition before calculating the efficiency
of the fixed local design.

We varied the prior standard deviation $\sigma$ from $0$ to $3$ in increments of $0.1$. We also included the values $\sigma=1/\sqrt{3}$ and $\sigma=5/\sqrt{3}$, corresponding to the slope variances $1/3$ and $100/12$, respectively. These variances are shared by the normal and uniform slope priors in Table~\ref{tab:priors_onefactor}. For the Poisson and logistic models, expectations over the slope were evaluated using Gauss-Hermite quadrature for normal priors and Gauss-Legendre quadrature for uniform priors, with $100$ nodes. The resulting designs were checked using $200$ nodes.
At $\sigma=0$, expectations were evaluated directly at $\beta_1=3$. The computed solutions satisfied $\max_x d_{\pi_\sigma}(x,\xi_{\mathrm B}^*)\leq2+2\times10^{-6}$.

\begin{figure}[H]
\centering
\includegraphics[width=\textwidth]{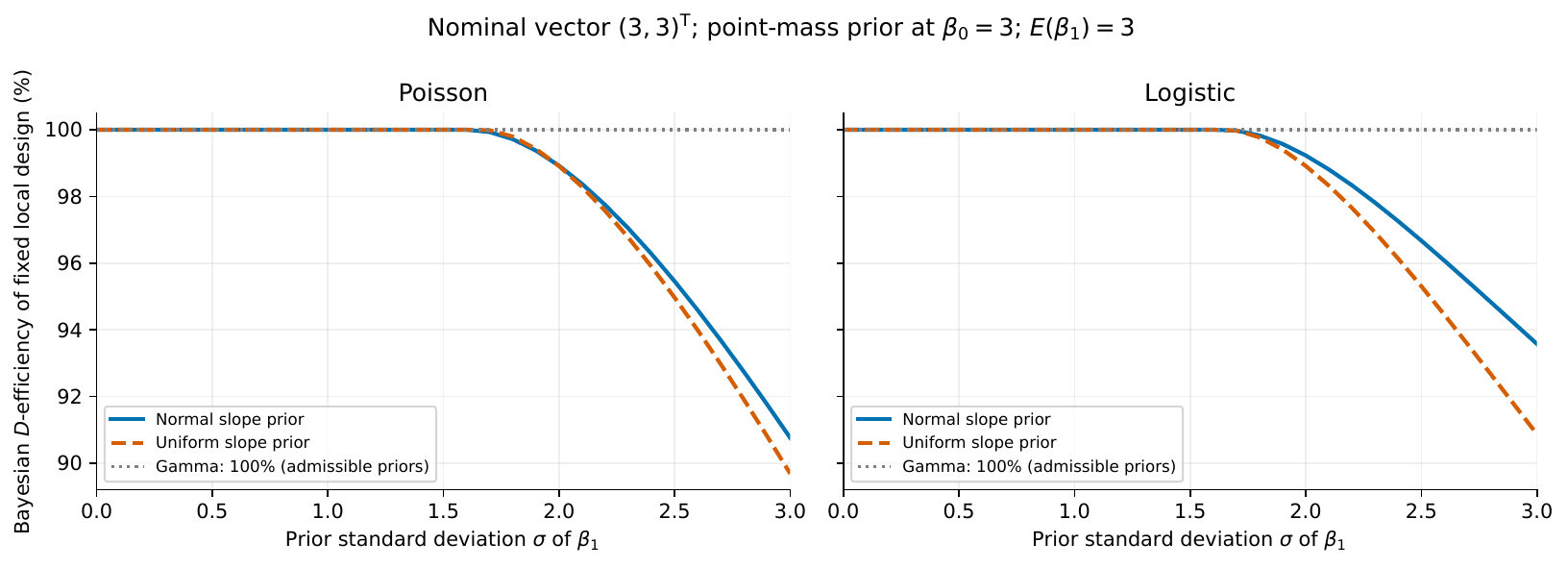}
\caption{Bayesian $D$-efficiency of the fixed locally $D$-optimal designs at $\widetilde{\BS{\beta}}=(3,3)^\trp$. The blue solid curve corresponds to $\beta_1\sim N(3,\sigma^2)$ and the orange dashed curve to $\beta_1\sim U(3-\sqrt{3}\sigma,3+\sqrt{3}\sigma)$, with a point-mass intercept prior at $3$. Each local design is compared with the Bayesian optimal design under the same prior. The dotted line gives the gamma benchmark under admissible priors.}
\label{fig:local_prior_efficiency}
\end{figure}

Figure~\ref{fig:local_prior_efficiency} shows the efficiency of the locally $D$-optimal design as prior dispersion increases. For Poisson regression under both normal and uniform slope priors, the local design remained Bayesian $D$-optimal over all designs at every tested value
$\sigma\leq1.6$, within the stated numerical tolerance. The first tested value yielding a three-point Bayesian
design was $\sigma=1.7$ for both prior families. For logistic model, the optimal support changes slightly for small values of $\sigma$, with negligible efficiency loss from retaining the local design. As $\sigma$ increases, this loss becomes
more substantial. The normal and uniform curves also show that priors with the same mean and variance can give different efficiencies. For example, at $\sigma=5/\sqrt{3}$, the efficiencies under the normal slope prior are $91.89\%$ for Poisson model and $94.29\%$ for logistic model. Under the uniform slope prior, they are $90.95\%$ and $91.89\%$, respectively. Thus, both the amount of uncertainty and the distributional form can affect the benefit of replacing the local design. At $\sigma=3$, the efficiency of the local Poisson design under the uniform slope prior is approximately $89.68\%$, corresponding to an efficiency loss of $10.32\%$.

\section{Discussion}
\label{sec5}

This work examined the robustness of locally $D$-optimal
designs under parameter uncertainty for single-factor
gamma, Poisson and logistic regression models.
For gamma regression under the links considered, the
equally weighted endpoint design is locally $D$-optimal
throughout the admissible parameter space and therefore
remains Bayesian $D$-optimal for every admissible prior.

The Poisson model showed different behavior. Within the
saturated class, local optimality is preserved when the
nominal parameter equals the prior mean. In this case,
the local and Bayesian $D$-optimal designs coincide.
When optimization is extended to all designs, however,
the numerical results showed departures from the local
design under some highly dispersed priors.

For logistic regression, the Bayesian optimal support
was sensitive to the form and dispersion of the prior.
Under concentrated priors, the changes in support were
small, and retaining the local design resulted in only
small efficiency losses. Increasing slope uncertainty
produced more noticeable changes in the optimal design
and reduced the efficiency of retaining the local design.

The results distinguish preservation of optimality from
retention of high Bayesian $D$-efficiency.
For Poisson and logistic regression, the numerical
comparisons concern the nominal parameter vector
$\widetilde{\BS{\beta}}=(3,3)^\trp$ and the prior
specifications considered. Within these scenarios,
the local designs retained high efficiency under
concentrated priors, whereas more dispersed priors
produced larger losses. In the slope-uncertainty study
with a point-mass intercept prior at 3, the local
Poisson design had an efficiency loss of $10.32\%$
under the uniform slope prior with $\sigma=3$.
These quantitative findings are specific to the
nominal parameter values and priors studied.

\section*{Declaration of Generative AI Use}
OpenAI's GPT-6 Astra was used under the author's direction to improve the readability and language of the manuscript. The authors take full responsibility for the content of the publication.

\section*{Conflict of Interest}
The  author states that there is no conflict of interest.

\bibliographystyle{elsarticle-harv}
\bibliography{reference}

\end{document}